\documentclass[preprint]{imsart}

\RequirePackage[OT1]{fontenc}
\RequirePackage{amsthm,amsmath, amssymb, enumerate}
\RequirePackage[numbers]{natbib}
\RequirePackage[colorlinks,citecolor=blue,urlcolor=blue]{hyperref}
\RequirePackage{hypernat}

\arxiv{math.PR/0000000}

\startlocaldefs
\numberwithin{equation}{section}
\theoremstyle{plain}
\newtheorem{theorem}{Theorem}[section]

\newtheorem{lemma}{Lemma}[section]
\newtheorem{definition}{Definition}[section]

\endlocaldefs

\newcommand{\R}{\mathbb R}
\newcommand{\N}{\mathbb N}
\newcommand{\Z}{\mathbb Z}

\begin{document}

\begin{frontmatter}
\title{Local Pinsker inequalities via
  Stein's discrete density approach}
\runtitle{Local Pinsker inequalities via
  Stein's discrete density approach}
\begin{aug}
\author{\fnms{Christophe}
  \snm{Ley}\thanksref[1]{1}\ead[label=e1]{chrisley@ulb.ac.be}} 
\and
\author{\fnms{Yvik}
  \snm{Swan}\ead[label=e2]{yvik.swan@uni.lu}}

\thankstext[1]{1}{Supported by a Mandat de Charg\'e de Recherche from
  the Fonds National de la Recherche Scientifique, Communaut\'e
  fran\c{c}aise de Belgique. Christophe Ley is also a member of ECARES.} 
\runauthor{C. Ley and Y. Swan}

\affiliation{E.C.A.R.E.S. and Universit\'e Libre de Bruxelles}
\address{
Universit\'e Libre de Bruxelles\\
D\'epartement de Math\'ematique\\
Boulevard du Triomphe \\
Campus Plaine -- CP210\\
B-1050 Brussels\\
\printead{e1}}
\address{
Universit\'e du Luxembourg\\
Facult\'e des Sciences, de la Technologie et de la Communication\\
 Unit\'e de Recherche en Math\'ematiques\\
6, rue Richard Coudenhove-Kalergi\\
L-1359 Luxembourg\\
\printead{e2}}

%
\end{aug}

\begin{abstract} Pinsker's inequality states that the relative entropy
  $d_{\mathrm{KL}}(X, Y)$ between two random variables $X$ and $Y$
  dominates the square of the total variation distance
  $d_{\mathrm{TV}}(X,Y)$ between $X$ and $Y$. In this paper we
  introduce generalized Fisher information distances $\mathcal{J}(X,
  Y)$ between discrete distributions $X$ and $Y$ and prove that these
  also dominate the square of the total variation distance. To this
  end we introduce a general discrete Stein operator for which we
  prove a useful covariance identity.  We illustrate our approach with
  several examples. Whenever competitor inequalities are available in the
  literature, the constants  in ours are at least as good, and, in
  several cases, better.
\end{abstract}



\end{frontmatter}

\section{Introduction}
\label{sec:stein-type-ident}

Let $X\sim p$ and $Y \sim q$ be two real-valued  random variables.  The
 \emph{relative entropy}  between $X$ and $Y$
 (a.k.a. \emph{Kullback-Leibler divergence}, see \cite{KullIEEE}) is
 defined as 
 \begin{equation}\label{eq:25}
   d_{\rm KL}(Y|| X) =  d_{\rm KL}(q|| p) = {\rm E}_q \left[ \log \left( \frac{q(Y)}{p(Y)} \right) \right],
 \end{equation}
 where $\mathrm{E}_q \left[ h(Y) \right]$ stands for the expectation
 of $h$ under $q$. Although not a \emph{bona fide} probability
distance  (absence of symmetry, no triangular inequality),  Gibbs'
inequality (see, e.g., \cite{CoTh06}) 
\begin{equation*}
  d_{\rm KL}(Y|| X) \ge 0 \mbox{ with equality if and only if } p=q
\end{equation*}
entails that $d_{KL}(Y || X)$ does indeed  quantify a particular form of
discrepancy (in terms of the entropies) between the law of $X$ and
that of $Y$.  Moreover, letting $d_{\rm TV} (X, Y)$ stand for the
\emph{total variation distance} between $p$ and $q$ (a precise
definition is given in Section~\ref{sec:inform-ineq-1}), the
\emph{Pinsker's inequality} (see, e.g., \cite{GS02,CoTh06}) 
\begin{equation}
  \label{eq:4}
  2 d_{\rm TV} (X, Y) \le   \sqrt{ 2   d_{\rm KL}(Y|| X)}
\end{equation}
implies  that  the relative entropy   dominates the total variation
distance and thus, also, a large class of classical
probability distances (including  the Wasserstein distance, see
e.g. \cite{GS02} for an overview of the interrelations between
probability metrics).

Fix $X=N$ a standard Gaussian random variable and consider absolutely
continuous random variables $Y$ with differentiable density $q$ and
finite variance which we set to 1. 
Estimates 
on $d_{\mathrm{KL}}(Y ||N) $ are typically obtained through control of
the \emph{Fisher information distance} (FID) between the law of $Y$
and the Gaussian, which is defined as
\begin{equation}\label{eq:24}
  J_N(Y) = \mathrm{E}_{q}\left[ \left( \frac{q'(Y)}{q(Y)}+Y \right)^2
  \right] = I(Y)-1,
\end{equation}
with $I(Y)= \mathrm{E}_{q}\left[ \left( {q'(Y)}/{q(Y)}\right)^2
  \right]$ the Fisher information of $Y$. The FID can be viewed as a ``local'' 
  version of the relative entropy (see, e.g.,
  \cite{BaBaNa03,MR2128239,carlen1991entropy,Br82,BA86}). Trivially
  positive, it satisfies
\begin{equation*}
J_N(Y) = 0 \mbox{ if and only if } Y
\stackrel{\mathcal{L}}{=} N
\end{equation*}
 so that $J_N(Y)$ indeed quantifies  discrepancy (this time in terms of
the Fisher informations) between $q$ and the Gaussian
distribution. Finally the FID
 dominates the total 
variation distance 
\begin{equation}\label{eq:shimi} 
  d_{\mathrm{TV}}(N,Y) \le  \sqrt{2 J_N(Y)}
\end{equation}
 (see \cite{Sh75, MR2128239})  so that (similarly as the relative
 entropy) proximity between the law of $Y$ and the Gaussian in terms of the Fisher information distance
implies proximity in terms of a wide variety of more classical
probability distances.

 Fix $X = Po(\lambda)$ a  rate-$\lambda$ Poisson random variable and consider
 discrete random variables $Y$ with probability mass function $q$ on the positive integers. There exist at least
 two ``local'' versions of \eqref{eq:25} which have been put to use in
 the literature on Poisson convergence, namely the \emph{discrete
   Fisher information}
\begin{equation}
  \label{generalized}
  \mathcal{J}(Po(\lambda), Y) := \mathrm{E}_q \left[ \left( \frac{\lambda
      q(Y-1)}{q(Y)} - Y \right)^2 \right]
\end{equation}
introduced in \cite{BaJoKoMa10} (itself a generalization of an information
functional presented in \cite{JM87}) and the \emph{scaled Fisher
  information}
\begin{equation}\label{scaled} 
\mathcal K(Po(\lambda),Y) := \lambda {\rm
  E}_q\left[\left(\frac{(Y+1)q(Y+1)}{\lambda q(Y)} -
    1\right)^2\right]
    \end{equation}
    introduced in \cite{KoHaJo05}. Both \eqref{generalized} and
    \eqref{scaled} are  trivially positive and 
    \[ \mathcal{J}(Po(\lambda), Y) = \mathcal K(Po(\lambda),Y) = 0
    \mbox{ if and only if } Y \sim Po(\lambda)\] so that these
    pseudo-distances indeed quantify a specific form of discrepancy between
    the density $q$ and the Poisson distribution. The scaled Fisher
    information $\mathcal{K}(Po(\lambda),Y)$ dominates the relative
    entropy $d_{\mathrm{KL}}(Po(\lambda),Y)$ (see \cite{KoHaJo05}) and thus, by Pinsker's
    inequality \eqref{eq:4}, 
    \begin{equation}\label{eq:konto} 
      d_{\mathrm{TV}}(Po(\lambda), Y) \le \sqrt{2 \mathcal{K}(Po(\lambda),Y)}.
    \end{equation}
Consequently, as above,  proximity in terms of the functional
$\mathcal{K}(Po(\lambda),Y)$ entails proximity in  terms of a wide variety of more classical
probability distances.

Inequalities \eqref{eq:shimi} and \eqref{eq:konto} are
local versions of  inequality \eqref{eq:4} with respect to a fixed
target distribution $X$. Moreover the three functionals \eqref{eq:24},
\eqref{generalized} and \eqref{scaled} are of the form 
\begin{equation*}
  \mathcal{J}(X, Y) = \mathrm{E}_q \left[ (r(p,q)(Y))^2 \right]
\end{equation*}
for $r(p,q)$ a mean-0 functional which we interpret as a score
function.  In view of the fact that Pinsker's inequality is valid
irrespective of the laws of $X$ and $Y$, it is natural to enquire
whether there exists some universal score function $r(p,q)$ whose
variance $\mathcal{J}(X,Y)$ provides an informative ``information
distance'' between the laws of $X$ and $Y$ such that (i)
$\mathcal{J}(X,Y)\ge0$ with equality if and only if $X
\stackrel{\mathcal{L}}{=}Y$, and (ii) $\mathcal{J}(X, Y)$ satisfies the local Pinsker's
inequality
\begin{equation}\label{eq:pinskoloc} 
  d_{\mathrm{TV}}(X, Y) \le \kappa \sqrt{\mathcal{J}(X,Y)}
\end{equation}
for $\kappa$ some constant whose value only depends on the properties
of the target distribution $p$.

A partial answer to this question is already known in case $p$ and $q$
are continuously differentiable probability density functions. Indeed  in
\cite{LS12a} we introduce the  \emph{generalized Fisher
information distance} 
\[
  \mathcal{J}(X, Y) := \mathrm{E}_q\left[ \left(
      \frac{p'(Y)}{p(Y)}-\frac{q'(Y)}{q(Y)} \right)^2 \right]
\]
which is a generalization of \eqref{eq:24} to arbitrary densities $p$
and $q$ (note how, if $p$ is the standard Gaussian density, we have
$p'(x)/p(x) = -x$ so that we recover $\mathcal{J}(N,Y) = J_N(Y)$ the FID). Under
assumptions on the supports of $p$ and $q$ we prove that $ \mathcal{J}(X, Y)$ 
satisfies \eqref{eq:pinskoloc}  and, for $p$ the Gaussian, 
recover the constant $\kappa_p = \sqrt2$, and thus inequality \eqref{eq:shimi}.

The approach developed in \cite{LS12a} is reserved to continuously
differentiable distributions on the real line, and the purpose of the
present paper is to cover the case of discrete distributions. Before
delving into the specifics of the discrete case, we start with an
intuitive overview of our approach.

\subsection{Sketch of the approach}
\label{sec:overview-approach}

 Fix $[a, b]=\left\{ a, a+1, \ldots, b \right\}$ a collection
of consecutive integers and
consider a random variable   $X \sim p$ with $p$ a discrete probability
distribution on $[a, b]$. Let $\Delta^{\eta}$ be the classical forward
($\eta=1$) or backward ($\eta=-1$) difference
operator on $\Z$  (see \eqref{eq:12} for a precise definition) 
and   define  the operator
$\mathcal{T}_p^{\eta}$  via the duality relationship 
\begin{equation}
  \label{eq:7}
  \mathrm{E}_p \left[  (\mathcal{T}_p^{\eta} f)(X) g(X) \right] =    \mathrm{E}_p \left[
    f(X) \Delta^{\eta} g(X) \right]
\end{equation}
which we require to hold for  all bounded functions  $g$ on $\Z$ and all   $f$
belonging to some class 
$\mathcal{F}(p)$ which satisfy 
the appropriate boundary conditions (see Definition \ref{def1}).
Setting
$g(x) = 1$ in \eqref{eq:7} we immediately deduce that
\begin{equation}\label{eq:expzero} 
   \mathrm{E}_p \left[ ( \mathcal{T}_p^{\eta} f)(X) \right]=0
\end{equation}
for all $f \in \mathcal{F}(p)$;  in Theorem \ref{theo:dda} we prove
that the converse also holds true, i.e. if $Y\sim q$ and $
\mathrm{E}_q \left[ ( \mathcal{T}_p^{\eta} f)(Y) \right]=0 $ for all
$f \in \mathcal{F}(p)$ then $p=q$.

  Operator $\mathcal{T}_p^{\eta}$ is a generalization of the so-called
\emph{Stein operators} from the literature on Stein's method
\cite{BaHoJa92,BC05,NP11,ChGoSh11} and the resulting characterization
(Theorem~\ref{theo:dda}) is a generalization of the so-called
\emph{density approach} adapted to the discrete setting, see
e.g. \cite{S01,GoRe12}. In the Appendix \ref{sec:examples}, we will
discuss specific examples for various choices of $p$ and show how our
operators contain many of the {Stein operators} which arise
through other (sometimes more complex) methods, see
e.g. \cite{Ho04}. 

The connection between Stein's method and information theory is
implicit in the works \cite{Jo04,BaJoKoMa10,Sh75} and is central to
\cite{LS12a,nourdin2013entropy}.  See also the works
\cite{sason2013improved,sason2012information,sason2012entropy,sason2012entropybern}
for alternative general considerations on the connexions between the
two topics in the discrete setting. In this work as well we make use
of a variation of this method, as follows.  Given $X\sim p$ and $Y
\sim q$ two random variables and $l$ some test function, consider the
solution $f_l^p$ of the difference equation (a.k.a. Stein equation)
\begin{equation}\label{eq:steak} 
  (\mathcal{T}_p^{\eta}f_l^p) (x) = l(x) - \mathrm{E}_p \left[ l(X) \right].
\end{equation}
Much is known, from the literature on Stein's method, on the
properties of the function $f_l^p$ for several choices of target $p$ (see, e.g., \cite{BaHoJa92}). 
Taking expectations (w.r.t. $q$) on both sides of \eqref{eq:steak}
and using fact  \eqref{eq:expzero} we get 
 \begin{align}
  {\rm E}_q[l(Y)] -{\rm E}_p[l(X)]&  =   {\rm E}_q [(
  \mathcal{T}_p^{\eta} f_l^p)(Y) ] \nonumber \\
  & = {\rm E}_q [(
  \mathcal{T}_p^{\eta}f_l^p)(Y) - (\mathcal{T}_q^{\eta}f_l^p)(Y) ]\label{eq:struu} 
\end{align}
under the assumption that $f_l^p \in \mathcal{F}(q)$. Furthermore, it is easy to prove
(see \eqref{eq:fac_prop_1})
that we have the decomposition
\begin{equation}
  \label{eq:decompooo}
   (\mathcal{T}_p^{\eta}f)(x) - (\mathcal{T}_q^{\eta}f)(x) = f(x) r^{\eta}(p, q)(x) + \epsilon(x)
\end{equation}
where $\epsilon$ has $q$-mean 0  and $r^{\eta}(p,q)$ is
some functional of the densities $p$ and $q$ (and not of
$f$) which, as we shall see, turns out to be a score function. Plugging \eqref{eq:decompooo} into \eqref{eq:struu} we 
get
\begin{equation}
  \label{eq:impo}
   {\rm E}_q[l(Y)] -{\rm E}_p[l(X)] = \mathrm{E}_q \left[ f_l^p(Y) r^{\eta}(p, q)(Y) \right].
\end{equation}
Now, many probability distances (total variation distance, Kolmogorov
distance, Wasserstein distance,...) can be written under the form
\begin{equation*}
  d_{\mathcal{H}}(X, Y) = \sup_{l\in \mathcal{H}} \left| {\rm E}_q[l(Y)] -{\rm E}_p[l(X)]  \right|
\end{equation*}
for $\mathcal{H}$ some class of functions (see, e.g., \cite[Appendix
C]{NP11}).  
Taking suprema on either side
of \eqref{eq:impo} we obtain
\begin{equation}
  \label{eq:16}
 d_{\mathcal{H}}(X, Y)  = \sup_{l\in \mathcal{H}}\mathrm{E}_q\left[
   \left| f_l^p(Y) r^{\eta}(p, q)(Y) \right|\right].
\end{equation}

We will use (Section \ref{sec:inform-ineq-1}) equality \eqref{eq:16}
to derive generalized Fisher information distances (for arbitrary
discrete distributions)  which we will  prove to satisfy  the local Pinsker's inequality
\eqref{eq:pinskoloc}  with an explicit constant
$\kappa$. In particular  we will introduce 
(i)  the discrete Fisher information distance
\begin{equation*}
  \mathcal{J}_{\mathrm{gen}}(X,Y) = \mathrm{E}_q \left[ \left( \frac{q(Y-1)}{q(Y)}-\frac{p(Y-1)}{p(Y)} \right)^2 \right]
\end{equation*}
(Section
\ref{sec:ineq-backw-diff}) which generalizes \eqref{generalized} and
(ii) the  scaled Fisher information distance 
\begin{equation*}
  \mathcal{K}_{\mathrm{gen}}(X, Y) = \mathrm{E}_q \left[ \left(
      \frac{p(Y)q(Y+1)}{p(Y+1)q(Y)}-1 \right)^2 \right]
\end{equation*}
 (Section
\ref{sec:ineq-forw-diff}) which generalizes \eqref{scaled}. These are
not the only discrete information distances that can be obtained by
our approach, although they are the most relevant in view of the
current literature on the topic. We illustrate (Section
\ref{sec:other-inequalities}) an alternative
construction in a  specific setting related to the recent reference
\cite{fulman2012stein}, and show that here as well our inequalities
are competitive.

\subsection{Outline of the paper}
\label{sec:outline-paper}


We start, in Section \ref{sec:results1}, by rigorously defining all
the concepts appearing in Section \ref{sec:overview-approach}. We also
provide explicit conditions under which the manipulations are
permitted. In Section \ref{sec:inform-ineq-1} we discuss the local
Pinsker's inequalities obtainable from \eqref{eq:16} and provide
several examples; we also compare our bounds with those already
available in the literature. Finally the Appendix contains details, proofs and
examples from Section \ref{sec:results1}.

\section{Stein's density approach for discrete distributions} \label{sec:results1}

Let $\mathcal{G}$ be the collection of probability mass functions 
$p:\Z\rightarrow[0,1]$ with  support $S_p:=\{x\in\Z:
p(x)>0\}$  a discrete interval $[a, b] := \{a, a+1,
\ldots, b\}$ for $a< b \in \Z\cup\{\pm \infty\}$. 
We will, 
   in the sequel,   abuse  language by referring to probability mass functions as
   (discrete) densities. Throughout we
adopt the convention that sums running over empty sets equal $0$, and
that
\begin{equation}
  \label{eq:10}
  \frac{1}{p(x)}=\left\{
\begin{array}{cl}
\frac{1}{p(x)}&\mbox{if}\,x\in S_p\\
0&\mbox{otherwise}.
\end{array}
\right.
\end{equation}
Note how, in particular, convention~\eqref{eq:10} implies that
$p(x)/p(x)=\mathbb{I}_{S_p}(x)$, the 
indicator of the support $S_p$.   
We will write
${\rm E}_p[l(X)]=\sum_{x\in S_p}l(x)p(x)$ for $p\in\mathcal{G}$
and $l$ a $p$-summable function. Furthermore we introduce 
the \emph{$\eta$-difference 
     operator} 
   \begin{equation}
     \label{eq:12}
     \Delta^{\eta}h(x) = \frac{1}{\eta} \left( h(x+\eta) - h(x) \right)
   \end{equation}
   for all functions $h$ taking their values on $\Z$. (Operators of
   the form  \eqref{eq:12} are not the only choice of ``discrete
   derivative operator''; see e.g. 
   \cite{hillion2011natural} for an alternative).   

\begin{definition}\label{def1}
Let $p\in\mathcal{G}$  and $|\eta |=1$. We define  (i) the collection
  $\mathcal{F}^{\eta}(p)$ of functions $f:\Z\rightarrow\R$ such that
$     \sum_{j=a}^b
  \Delta^{\eta}(f(j)p(j))=0,
$
  and (ii) the operator
  $\mathcal{T}_p^{\eta}:\mathcal{F}^{\eta}(p)\rightarrow\Z^*:f\mapsto
  \mathcal{T}_p^{\eta}f$ given by
\begin{equation}\label{operator}
\mathcal{T}_p^\eta f:\Z\rightarrow\R:x\mapsto \mathcal{T}_p^\eta
f(x):=\frac{1}{p(x)}\Delta^\eta(f(x)p(x)). 
\end{equation} 
We call $\mathcal{F}^{\eta}(p)$ the class of $\eta$-\emph{test functions} associated
with $p$, and $\mathcal{T}_p^{\eta}$ the $\eta$-\emph{Stein operator} associated
with $p$.
\end{definition}

The first condition in  Definition \ref{def1}  (control of the
functions at the edges of the support) ensures that we have the
integration by parts formula 
\begin{equation*}
\mathrm{E}_p \left[ (\mathcal{T}_p^{\eta}f)(X) g(X) \right]  = -
\mathrm{E}_p \left[ f(X)\Delta^{\eta}g(X)  \right] 
\end{equation*}
for all functions $g$ for which the above makes sense. 

In particular, the class $\mathcal{F}^\eta(p)$ is tailored to ensure  that
$ {\rm E}_p[\mathcal{T}_p^{\eta}f(X)]=0$ for all $f \in
\mathcal{F}^\eta(p)$.  The following result (whose proof is deferred to 
the Appendix) shows that the converse holds true as well. 
 
\begin{theorem}[Discrete density approach] \label{theo:dda} Fix
  $|\eta|=1$ and let $X$ be
  a discrete random variable with density $p\in\mathcal{G}$. Let $Y$
  be another discrete random variable with density $q\in\mathcal{G}$. Then ${\rm E}_q[\mathcal{T}_p^{\eta}f(Y)]=0$ for all
  $f\in\mathcal{F}^\eta(p)$ if, and only if, either ${\rm P}(Y\in S_p)=0$
  or ${\rm P}(Y\in S_p)>0$ and ${\rm P}(Y\leq z\,|\,Y\in S_p)={\rm
    P}(X\leq z)$ for all $z\in S_p$.
\end{theorem}
 
Theorem \ref{theo:dda} is  a general \emph{Stein
  characterization}. 
Expounding, for $\eta=1$, the forward difference in
\eqref{operator} we get the same expression as \cite[Equation
(8)]{GoRe12}. Our density approach and theirs are not equivalent, as
described in \cite[Remark 2.1]{GoRe12}.  The differences between their
assumptions and ours are due to the ``difference of a product''
structure of \eqref{operator}.  Examples wherein we apply Theorem
\ref{theo:dda} to specific choices of $p$ and further details are
discussed in the Appendix.

Fix, for the sake of convenience, $S_p = [0, \ldots, M]$ and
$S_q=[0,\ldots,N]$, for some integers $0 \le N \le M \le \infty$. Note
in particular that we hereby ensure the crucial assumption $S_q
\subseteq S_p$. Now suppose that $\mathcal{F}^{\eta}(p) \cap
\mathcal{F}^{\eta}(q) \neq \emptyset$ and choose some $f$ in this
intersection. Then, for this $f$, we can write
\begin{align}
\label{eq:fac_prop_1} 
&  \mathcal{T}_p^{\eta}f(x)  = \mathcal{T}_q^{\eta}f(x) +
  \mathcal{T}_p^{\eta}f(x) - \mathcal{T}_q^{\eta}f(x) \nonumber \\
& = \mathcal{T}_q^{\eta}f(x) + \frac{1}{\eta}\left(
  \frac{\Delta^{\eta}(f(x)p(x)) }{p(x)} -
  \frac{\Delta^{\eta}(f(x)q(x)) }{q(x)}\right)\nonumber \\
& = \mathcal{T}_q^{\eta}f(x) + f(x+\eta)\frac{1}{\eta}\left(
  \frac{p(x+\eta)}{p(x)}-  \frac{q(x+\eta)}{q(x)} \right)-\frac{1}{\eta}
  f(x) \mathbb{I}_{[N+1,  \ldots, M]}(x),   
\end{align}
where the indicator function equals 0 if $M=N$.
Next let $l:\Z\to\R$ be a function such that both ${\rm E}_p[l(X)]$ and
${\rm E}_q[l(Y)]$ exist and consider the solution $f_{l}^{p, \eta}$ of
the
{difference (Stein) equation}  
\begin{equation}
  \label{eq:stein}
  \mathcal{T}_p^{\eta}f(x) = l(x) - {\rm E}_p[l(X)]. 
\end{equation}
 As in the proof of Theorem \ref{theo:dda} (see identities \eqref{eq:20}
and \eqref{eq:21}) it is easy to show that 
the solutions to \eqref{eq:stein} are given by 
\begin{equation}\label{eq:steineq_sol3} 
f_l^{p, +}: \Z \to \R : x \mapsto \sum_{k=0}^{x-1}(l(k)-
  {\rm E}_p[l(X)]) \frac{p(k)}{p(x)}\end{equation}
for $\eta=1$ (the forward difference operator) and 
\begin{equation}\label{eq:steineq_sol4} 
f_l^{p, -}: \Z \to \R : x \mapsto \sum_{k=0}^{x}(l(k)-
  {\rm E}_p[l(X)]) \frac{p(k)}{p(x)}\end{equation}
for $\eta=-1$ (the backward difference operator).  Recall that empty
sums  are set to 0. The functions $f_l^{p, \eta}$ as 
defined above trivially belong to 
$\mathcal{F}^{\eta}(p)$.  

To pursue we  need the   following assumption.

\

\noindent \emph{Assumption~A :} The distributions $p$ and $q$ are such
that  the solutions $f_l^{p,
  \eta}$ of the Stein equation \eqref{eq:stein} satisfy $f_l^{p,
  \eta} \in \mathcal{F}^{\eta}(p) \cap \mathcal{F}^{\eta}(q)$ for $|\eta|=1$.

\

\noindent For any given target $p$ it is easy to determine conditions
on $q$ and $l$ for Assumption~A to be satisfied. These conditions are
not restrictive. 

Under Assumption~A we can take expectations with respect to $q$
on either sides of~(\ref{eq:fac_prop_1}) applied to a solution of
~(\ref{eq:stein}). Since  $S_q\,\cap\, [N+1,\ldots, M]=\emptyset$ we
have $ \mathbb{I}_{[N+1,\ldots, M]}(Y) = 0$. Also $  \mathrm{E}_q
\left[ \mathcal T_q^{\eta}f_l^{p, \eta}(Y) \right]=0
$,
through Theorem~\ref{theo:dda} since
$f_l^{p,\eta}\in\mathcal{F}^\eta(q)$ by Assumption~A. Hence 
\begin{align*}
 {\rm E}_q[l(Y)] - {\rm E}_p[l(X)]
&  = \mathrm{E}_q \left[ \mathcal T_q^{\eta}f_l^{p,\eta}(Y) \right]  +
\mathrm{E}_q \left[ f_l^{p, \eta}(Y+ \eta) r^{\eta}(p,q)(Y)  \right],
\end{align*}  
with
\begin{equation} \label{eq:gensssscor}  r^{\eta}(p,q)(x) :=
   \frac{1}{\eta}\left(\frac{p(x+\eta)}{p(x)} - \frac{q(x+\eta)}{q(x)}\right).\end{equation}
We have proved  the following result.

\begin{lemma}\label{prop:steinentropy2}   
Take $p, q \in \mathcal{G}$ with $S_q\subseteq S_p$ and $l:\Z\to\R$  a function such that
${\rm E}_p[l(X)]$ and ${\rm E}_q[l(Y)]$ exist.  Suppose moreover that
Assumption~A holds. Then 
\begin{equation} \label{eq:fundeq2}  {\rm E}_q[l(Y)] -{\rm E}_p[l(X)] =  
{\rm E}_q [f_{l}^{p, \eta}(Y+\eta) r^{\eta}(p,q)(Y)], 
\end{equation}
with $f_l^{p, \eta}$ as in \eqref{eq:steineq_sol3} and
\eqref{eq:steineq_sol4} and $ r^{\eta}(p,q)$ as in \eqref{eq:gensssscor}.
 \end{lemma}
Following the terminology from
 \cite{GoRe12,APP11,APP07} we call \eqref{eq:fundeq2} a Stein (or Stein-type)
 identity.
Similarly as its counterpart \cite[Lemma 3.2]{LS12a} in the absolutely continuous setting,
 Lemma \ref{prop:steinentropy2}  provides the
 connection between our version of the discrete density approach
 from Theorem \ref{theo:dda} and discrete
 information inequalities.

\section{Local   Pinsker inequalities}
\label{sec:inform-ineq-1}

As already mentioned in the introduction, a wide variety of probability metrics can be written under the form
\begin{equation}\label{eq:hdizst}d_{\mathcal H} (X, Y) = \sup_{l \in
    \mathcal{H}} \left|{\rm E}_q[l(Y)] - {\rm E}_p[l(X)]
  \right|\end{equation} 
for some class of functions $\mathcal{H}$. 
In particular the total variation distance
\begin{equation*}
  d_{\rm TV}(X,Y):=\frac{1}{2}\sum_{x\in\N}\left|p (x)-q(x)\right| =  \sup_{l\in\{h\}}\left|{\rm E}_q \left[l(Y)\right] -
  {\rm E}_{p}  \left[l(X)\right]\right|,
\end{equation*}
where the supremum in the second equality  is taken over a set
containing one single function, namely   
\begin{equation*}
  h(x) := 
\frac{1}{2} \left( \mathbb{I}_{[p(x)\le q(x)]} -
  \mathbb{I}_{[p(x)\ge q(x)]} \right) =
\mathbb{I}_{[p(x)\le q(x)]} - \frac{1}{2}.
\end{equation*}
 Other distances
such as the Kolmogorov, the Wasserstein, the supremum-distance or the
$L^1$-distance can also be written under <  the form \eqref{eq:hdizst} -- we refer the reader to
\cite{GS02} or to  \cite[Appendix C]{NP11} for an overview. 

In view of \eqref{eq:hdizst}, it
is    natural to take  suprema on either side
of  \eqref{eq:fundeq2}   to deduce  that, whenever Assumption~A is
satisfied, we have
\begin{equation}
  \label{eq:11}
  d_{\mathcal H} (X, Y) =  \sup_{l \in
  \mathcal{H}} \left|\mathrm{E}_q \left[  f_l^{p, \eta}(Y+ \eta)
    r^{\eta}(p,q)(Y)  \right]\right|.
\end{equation}
Equation \eqref{eq:11} is a very powerful identity as it  
permits to identify \emph{natural}
discrete information distances which uniformly dominate all
probability distances of the form \eqref{eq:hdizst} through an
inequality in which only the constant is
distance-dependent. These  inequalities being 
valid for virtually any choice $(p, q)$,  we contend that their scope
is comparable with that of Pinsker's inequality \eqref{eq:4}, this
time for local versions of the (discrete) Kullback-Leibler divergence
\eqref{eq:25}.

\subsection{Fisher information inequalities via the backward difference operator}
\label{sec:ineq-backw-diff}

Choose the backward difference operator obtained for $\eta = -1$. Identity \eqref{eq:fundeq2} spells out as
\begin{equation}
  \label{eq:15}
   {\rm E}_q[l(Y)] -{\rm E}_p[l(X)] =  
{\rm E}_q [f_{l}^{p,-}(Y-1) r^{-}(p,q)(Y)]
\end{equation}
with $  r^{-}(p,q)(x) =\frac{q(x-1)}{q(x)}-\frac{p(x-1)}{p(x)}$ and 
with $f_{l}^{p,-}$ as in \eqref{eq:steineq_sol4}. 
 Taking suprema on either
side of \eqref{eq:15} and applying Cauchy-Schwarz we obtain the
following. 
\begin{theorem}\label{the:inform-ineq-discr}
  Take $p, q \in \mathcal{G}$ with $S_q\subseteq S_p$ and such that 
  $\mathcal{F}^{-}(p) \cap \mathcal{F}^{-}(q) \neq \emptyset$.
Let $d_{\mathcal H} (X, Y)$ be defined as in \eqref{eq:hdizst} for
some class of functions $\mathcal{H}$, and suppose that for all $l\in
\mathcal{H}$ the function $f_l^{p,-}$  
defined in \eqref{eq:steineq_sol4} exists and satisfies $f_l^{p,-} \in
\mathcal{F}^{-}(p) \cap \mathcal{F}^{-}(q)$. Then 
\begin{equation*}
\label{eq:29}
  d_{\mathcal H} (X, Y) \le\kappa_{\mathcal{H}}^{p, -} \sqrt{\mathcal J_{\rm gen}(X,Y)},
\end{equation*}
where 
\begin{equation}
  \label{eq:31}
  \mathcal J_{\rm gen}(X,Y) := {\rm
  E}_q\left[\left(\frac{q(Y-1)}{q(Y)}-\frac{p(Y-1)}{p(Y)}\right)^2\right]
\end{equation}
is the \emph{generalized discrete  Fisher information distance} between
the densities $p$ and $q$, and 
\begin{equation*}\label{eq:18}
 \kappa_{\mathcal{H}}^{p, -} :=  \sup_{l \in \mathcal{H}}\sqrt{{\rm
     E}_q \left[\left(f_l^{p,-}(Y-1)\right)^2\right]}. 
\end{equation*}
\end{theorem}

As an application suppose that $p$ and $q$ share the same
support.  Then we can write 
\begin{align*}
   \frac{q(x-1)}{q(x)} - \frac{p(x-1)}{p(x)} & =
   \frac{\Delta^{-}p(x)}{p(x)} - \frac{\Delta^{-}q(x)}{q(x)}
\end{align*}
so that  \eqref{eq:31} becomes 
\begin{equation}
  \label{eq:17}
  \mathcal J_{\rm gen}(X,Y) = {\rm
  E}_q\left[\left( \frac{\Delta^{-}p(Y)}{p(Y)} -
    \frac{\Delta^{-}q(Y)}{q(Y)}\right)^2\right]. 
\end{equation}
 The distance  \eqref{eq:17} extends  the Fisher information
 distance \eqref{generalized} to the comparison of any pair of
 densities $p,q$. Taking, in particular, $p$ a Poisson target we retrieve
\begin{align*}
  \mathcal{J}_{\mathrm{gen}}(Po(\lambda), Y) & = {\rm
  E}_q\left[\left( \left( 1- \frac{Y}{\lambda} \right)- 
    \frac{\Delta^{-}q(Y)}{q(Y)}\right)^2\right]\\
&= \frac{1}{\lambda^2}{\rm
  E}_q\left[\left( Y - 
    \frac{\lambda q(Y-1)}{q(Y)}\right)^2\right],
 \end{align*}
which in turn can be expressed as $ \frac{\sigma^2}{\lambda^2}-\frac{2}{\lambda}+ I(Y)$ 
with 
\begin{equation*}
  I(Y) = \mathrm{E}_q \left[ \left( \frac{\Delta^{-}q(Y)}{q(Y)}
  \right)^2 \right]
\end{equation*}
 the functional proposed in \cite{JM87} and
$\lambda, 
\sigma^2$ the mean and variance of $q$  
(see  also \cite[equation~3.1]{BaJoKoMa10}). In the
particular case of the 
Poisson distribution, the function $f_l^{p_{\lambda},-}(x-1)/\lambda$ is none other than the
usual solution of the standard equation  \eqref{eq:8} for which we
know  (see \cite[Theorem
2.3]{E05}) the estimate
\begin{align*}\label{eq:32}
  \left\|\frac{f_l^{p_{\lambda},-}(x-1)}{\lambda}\right\|_\infty \le \left( 1- \sqrt{\frac{2}{e\lambda }}\right)
\left(   \sup_{i \in \N} l(i) - \inf_{i \in \N} l(i)  \right);
\end{align*}
this is useful when $l$ is bounded as is the case, e.g., for the total
variation distance. Moreover, this boundedness of
$f_l^{p_{\lambda},-}$ also ensures that Assumption~A is satisfied
whatever $q$ (with support $\N$) we use, hence
Theorem~\ref{the:inform-ineq-discr} can be applied. Since we always
have
\begin{equation*}
   \kappa_{\mathcal{H}}^{p,-}  \le \sup_{l \in \mathcal{H}} \|f_l^{p,-}\|_{\infty},
\end{equation*}
we conclude from Theorem~\ref{the:inform-ineq-discr}  the information inequality 
\begin{equation*}
  \label{eq:36}
  d_{\rm TV}(Po(\lambda), Y) \le \left( 1- \sqrt{\frac{2}{e
        \lambda}}\right) \sqrt{\sigma^2-2\lambda+ \lambda^2I(Y)}.
\end{equation*}
Note that, for $q =p_{\lambda}$,  $I(Y) = \frac{1}{\lambda}$ and $\sigma^2=\lambda$ so that  $\sigma^2-2\lambda+
\lambda^2I(Y)=0$, as expected.

The information distance \eqref{eq:31} bears the defaults of its
originator \eqref{generalized} : if $p$ and $q$ do not share the same
support then $\mathcal{J}_{\mathrm{gen}}$ is infinite. In particular in the
Poisson case the quantity  for $q$ with
bounded support then, for some $k>0$, we have $q(k)
>0$ with $q(x+1) = 0$ so that $I(Y)=+\infty$ (see e.g. the discussion
at the beginning of  \cite[Section III]{KoHaJo05}). One way to avoid
this pathology is through a change in the derivative \eqref{eq:12}, as
follows. 

\subsection{Fisher information inequalities for the forward difference operator}
\label{sec:ineq-forw-diff}

Choose the forward difference operator, that is take \eqref{eq:12}
this time  with $\eta = 1$. Then  $r^+(p,q)(x)=\frac{p(x+1)}{p(x)}-\frac{q(x+1)}{q(x)}$
and $f_l^{p,+}$ is of the form~\eqref{eq:steineq_sol3}. If the target distribution
$p$ has support $\N$ then $p(x)/p(x+1)$ is finite  for all $x \in \N$ and the factorization
\begin{equation}\label{eq:factorrr} 
  f_l^{p,+}(x+1)r^+(p,q)(x) =  \left\{ f_l^{p,+}(x+1)\frac{p(x+1)}{p(x)} \right\}\left(1-\frac{q(x+1)p(x)}{q(x)p(x+1)}\right)
\end{equation}
is well-defined for all $x$. We introduce the scaled score function 
\begin{equation}
  \label{eq:scascorfun}
r_{\mathrm{sca}}(p, q)(x) =   1-\frac{q(x+1)p(x)}{q(x)p(x+1)}
\end{equation}
and the analog of Theorem~\ref{the:inform-ineq-discr} is obtained by yet
another simple application of the Cauchy-Schwarz inequality to this
factorization.

\begin{theorem}\label{disc2}
Take $p, q \in \mathcal{G}$ with $S_q\subseteq S_p$ and such that 
  $\mathcal{F}^{+}(p) \cap \mathcal{F}^{+}(q) \neq \emptyset$.
Let $d_{\mathcal H} (X, Y)$ be defined as in \eqref{eq:hdizst} for
some class of functions $\mathcal{H}$, and suppose that for all $l\in
\mathcal{H}$ the function $f_l^{p,+}$, as 
defined in \eqref{eq:steineq_sol3}, exists and satisfies $f_l^{p,+} \in
\mathcal{F}^{+}(p) \cap \mathcal{F}^{+}(q)$. Then 
\begin{equation}
  \label{disc2_1}
  d_{\mathcal H} (X, Y) \le\kappa_{\mathcal{H}}^{p,+} \sqrt{\mathcal K_{\rm gen}(X,Y)},
\end{equation}
where 
\begin{equation*}
 \kappa_{\mathcal{H}}^{p,+} :=  \sup_{l \in \mathcal{H}}\sqrt{{\rm E}_q \left[\left(f_l^{p,+}(Y+1)\frac{p(Y+1)}{p(Y)}\right)^2\right]}
\end{equation*}
and 
$$
\mathcal K_{\rm gen}(X,Y) = \mathrm{E}\left[ (r_{\mathrm{sca}}(p,
  q)(Y))^2 \right] =  {\rm
  E}_q\left[\left(\frac{p(Y)q(Y+1)}{p(Y+1) q(Y)} -
    1\right)^2\right]
    $$
is the \emph{generalized scaled Fisher information} between
the densities $p$ and $q$.   
\end{theorem}

In the case  $p = Po(\lambda)$ we have
$p_\lambda(x+1)/p_\lambda(x)=\lambda/(x+1)$ so that 
\eqref{disc2_1} becomes
 \begin{align*}
 d_{\mathcal H} (Po(\lambda), Y) &\le\sup_{l \in \mathcal{H}}\sqrt{{\rm
     E}_q
   \left[\left(f_l^{p_\lambda,+}(Y+1)\frac{\lambda}{Y+1}\right)^2\right]}
 \sqrt{\mathcal K_{\rm gen}(Po(\lambda),Y)} \label{ours2}\\
&=\sup_{l \in \mathcal{H}}\sqrt{{\rm
     E}_q
   \left[\left(f_l^{p_\lambda,+}(Y+1)\frac{\sqrt{\lambda}}{Y+1}\right)^2\right]}
 \sqrt{\mathcal K(Po(\lambda),Y)} \nonumber
\end{align*}
with $\mathcal K(Po(\lambda),Y)=\lambda \mathcal{K}_{\rm gen}(Po(\lambda),Y)$ the scaled Fisher information
distance \eqref{scaled}.  Using a Poincar\'e inequality, \cite{KoHaJo05} show that,
for $q$ a discrete 
distribution with mean $\lambda$,    
\begin{equation} 
  \label{eq:39}
  d_{\rm TV}(Po(\lambda),Y)\le \sqrt{2 \mathcal K(Po(\lambda),Y)}.
\end{equation}
Our Theorem~\ref{disc2}  allows to improve on this result, through the inequality  (see again \cite[Theorem
2.3]{E05})
\begin{align*}
  \left\|\frac{f_l^{p_{\lambda},+}(x+1)}{x+1}\right\|_\infty \le \left( 1- \sqrt{\frac{2}{e\lambda }}\right)
\left(   \sup_{i \in \N} l(i) - \inf_{i \in \N} l(i)  \right); 
\end{align*}
indeed, this inequality combined with Theorem~\ref{disc2} yields
(under the appropriate and more general conditions than in
\cite{KoHaJo05})
\begin{equation}\label{eq:33}
d_{\rm TV}(Po(\lambda),Y)\le \sqrt{\lambda}\left(1 \wedge
  \sqrt{\frac{2}{e {\lambda}}}\right) \sqrt{\mathcal
  K(Po(\lambda),Y)}.\end{equation}   
For $\lambda<2/e$, we get $1 \wedge \sqrt{\frac{2}{e {\lambda}}}=1$
and hence the constant in \eqref{eq:33} is
$\sqrt{\lambda}<\sqrt{2/e}$; in case $\lambda>2/e$, this constant
equals $\sqrt{2/e}$. In both cases our constants improve on those from
\eqref{eq:39}. More generally one easily sees that, for instance, in
all examples considered in \cite{KoHaJo05} our constants are better.

\subsection{Other inequalities}
\label{sec:other-inequalities}

In certain cases  it is better  to work
directly from the Stein identity \eqref{eq:11} without applying the
Cauchy-Schwarz inequality.  We illustrate this in
the specific case  of
approximation of the rank distribution of random matrices over finite
fields, as studied recently in \cite{fulman2012stein}.

Let $M_n$ be chosen uniformly from
 $Mat(n,\theta)$ the collection of all $n\times n$ matrices over the
 finite field $\mathbb{F}_{\theta}$ of size $\theta\ge 2$. Let
 $Q_{\theta}^n = n-\mbox{rank}(M_n)$ and let   $Q_{\theta}$  be its limiting version as $n
 \to \infty$. Both  the distribution of $Q_{\theta}^n$
 ($q_{k, n}, k = 0, \ldots, n$,  say) and that of $Q_{\theta}$   ($q_k, k\ge0$, say) are known -- see \cite[equations (1),
 (2)]{fulman2012stein}. 
These distributions satisfy the recurrence relations 
\begin{equation*}
  \frac{q_{k-1}}{q_k} = \frac{(\theta^k-1)^2}{\theta}, \, k \in \N \mbox{ and }
  \frac{q_{k-1, n}}{q_{k,n}} = \frac{(\theta^k-1)^2}{\theta(1-\theta^{-n+k-1})}, \, k
  \in [0,n].
\end{equation*}
Using \eqref{eq:11} with forward difference $\Delta^{+}$ and   factorization
\eqref{eq:factorrr},  the corresponding  score function 
\eqref{eq:scascorfun} simplifies to  (for $p = q_k$ and $q = q_{k,n}$)
\begin{equation*}
  r_{\mathrm{sca}}(Q_{\theta}, Q_{\theta}^n)(x) = \theta^{-n+x}
\end{equation*}
so that 
\begin{equation}
  \label{eq:goldful}
  d_{\mathcal{H}}(Q_{\theta}, Q_{\theta}^n) = \sup_{l \in
    \mathcal{H}} \left| \mathrm{E} \left[ f_l^{\theta,
        +}(Q_{\theta}^n+1) \frac{\theta}{(\theta^{Q_{\theta}^n+1}-1)^2}
      \theta^{-n+Q_{\theta}^n} \right] \right|
\end{equation}
with 
$f_l^{\theta, +}$ the solution to the
difference equation   \eqref{eq:stein}  given by
\eqref{eq:steineq_sol3}. See  Appendix \ref{sec:examples} where we
outline the setup of Stein's method via our Theorem \ref{theo:dda}
applied to this choice of distribution. 

  Inequality \eqref{eq:goldful} allows to
recover the upper bound from  \cite[Theorem
1.1]{fulman2012stein}. Indeed  it is shown there  \cite[Lemma 3.3]{fulman2012stein}  that 
\begin{equation*}
 \mathrm{E} \left[ \theta^{Q_{\theta}^n} \right] =
 2-\frac{1}{\theta^n} \mbox{ and }
\left\|\frac{f_l^{\theta,
    +}(x+1)}{(\theta^{x+1}-1)^2}\right\|_{\infty}
\le  \frac{1}{\theta^2}+\frac{1}{\theta^3}, 
\end{equation*}
 if $l$ is an indicator
function. 
Plugging these facts into  \eqref{eq:goldful} we get
\begin{align*}
 d_{\mathrm{TV}}(Q_{\theta}, Q_{\theta}^n) & \le  \left(
    \frac{1}{\theta^2}+ \frac{1}{\theta^3} \right)
  \theta^{-n+1}  \mathrm{E} \left[ \theta^{Q_{\theta}^n} \right] \\
  & = \left(
    \frac{1}{\theta^2}+ \frac{1}{\theta^3} \right)
  \theta^{-n+1}  \left( 2-\frac{1}{\theta^n} \right)  \\
& \le
  \frac{2(1+1/\theta)}{\theta^{n+1}} \le \frac{3}{\theta^{n+1}}
\end{align*}
for all $\theta\ge 2$; this is the upper bound from \cite[Theorem 1.1]{fulman2012stein}.

One can also, using H\"older's inequality in \eqref{eq:goldful},
obtain bounds on the total variation distance in terms of higher
moments $\mathrm{E} \left[ \theta^{k Q_{\theta}^n} \right]$,
$k\ge1$. Initial computations show that the resulting inequalities are
of equivalent rate but with constants depending  on $\theta$ and  bigger than
3. It would be interesting to enquire whether better inequalities are
obtainable by exploiting the flexibility in \eqref{eq:11}. This is
outside of the scope of the present article.

\appendix

\section{Details from Section \ref{sec:results1}}
\label{sec:exampl-stein-oper}

\subsection{Proof of Theorem \ref{theo:dda}}
\label{sec:proofs}
 If ${\rm P}(Y\in S_p)=0$, the equivalence holds
  trivially so that we can take ${\rm P}(Y\in S_p)>0$. 
  We first check sufficiency. The equality ${\rm P}(Y\leq z\,|\,Y\in
  S_p)={\rm P}(X\leq z)$ for all $z\in S_p$ can be rewritten as ${\rm
    P}(Y=z)={\rm P}(X=z){\rm P}(Y\in S_p)$, hence as $q(z)=p(z){\rm
    P}(Y\in S_p)$, for all $z\in S_p$. Bearing in mind that the
  operator $\mathcal{T}_p^{\eta}f(x)=0$ for all $x \notin S_p$, the
  sufficiency is easily established through
\begin{align*}
{\rm E}_q[\mathcal{T}_p^{\eta}f(Y)]&
={\rm P}(Y\in S_p)\sum_{x\in S_p}\Delta^\eta(f(x)p(x))=0,
\end{align*}
the last equality following by definition of the class $\mathcal{F}^\eta(p)$. 
Next, to see the necessity, 
define, for $z\in \Z$, the functions $l_z(k):= ({\mathbb{I}}_{(-\infty, z]\cap\Z}(k) - {\rm P}(X \le
z))\mathbb{I}_{S_p}(k)$ for $k \in \Z$ and define
\begin{equation}
  \label{eq:20}
  f_z^{p, +1}: \Z\to \R: x \mapsto \frac{1}{p(x)} \sum_{k=a}^{x-1} l_z(k)p(k)
\end{equation}
and 
\begin{equation}
  \label{eq:21}
    f_z^{p, -1}: \Z\to \R: x \mapsto \frac{1}{p(x)} \sum_{k=a}^{x} l_z(k)p(k).
\end{equation}
Clearly these functions satisfy $  \Delta^\eta(f_z^{p,  \eta}(x)p(x))= l_z(x)p(x)$
so that, in particular, $f_z^{p, \eta}\in\mathcal{F}^{\eta}(p)$ and   
\begin{equation*}
  \mathcal{T}_{p}^{\eta}f_z^{p, \eta}(x)=l_z(x)
\end{equation*}
 for all $x\in S_p$. 
 Consequently, for this choice of test function we obtain
\begin{align*}
  \sum_{x\in S_p}\mathcal{T}_p^\eta f_z^{p, \eta}(x)q(x) & =\sum_{x\in
    S_p}l_z(x)q(x)\\
& ={\rm P}(Y\leq z\cap Y\in S_p)-{\rm P}(Y\in S_p){\rm P}(X\leq z),
\end{align*}
which, in combination with the hypothesis ${\rm
  E}_q\left[\mathcal{T}_p^{\eta}f_z^{p, \eta}(Y)\right] =0$, finally
yields ${\rm P}(Y\leq z\,|\,Y\in S_p)={\rm P}(X\leq z)$ for all $z\in
S_p$, whence the claim.

\subsection{Examples of Stein operators}
\label{sec:examples}

Theorem \ref{theo:dda} extends and
    unifies many corresponding results from the literature, as will be
    shown through the following examples.

Take $p(x)=p_\lambda(x)
$ the density of a mean-$\lambda$ Poisson random variable. Then the
class $\mathcal{F}^{+}(p) =: \mathcal{F}^{+}(\lambda)$ is composed of
all functions $f:\Z\rightarrow\R$ such that (i) $x\mapsto
\Delta^{+}(f(x)p_\lambda(x))$ is summable over $\N$ and (ii)
$f(0)p_\lambda(0)=\lim_{x\rightarrow\infty}f(x)p_\lambda(x)$ (which in
most cases equals 0). In particular, $\mathcal{F}^{+}(\lambda)$
contains the set of bounded functions $f$ such that $f(0)=0$ (this border requirement is
necessary in order to belong to $\mathcal{F}^+$, see Definition \ref{def1}(i)), for
which simple computations show that
$$\mathcal T_{\lambda}^{+}f(x) = \left(\frac{\lambda}{x+1}f(x+1)
  - f(x) \right)\mathbb{I}_\N(x).$$ This operator coincides with that
discussed in \cite[page 6]{GoRe12}. One could also consider only
functions of the form $f(x) =xf_0(x)$ for $f_0$ such that $x \mapsto
xf_0(x) \in \mathcal{F}^{+}(\lambda)$ in which case no restriction on
$f_0(0)$ (other than that it be finite) is then necessary to ensure
the required border behaviour. Plugging such functions into
\eqref{operator} and simplifying accordingly we obtain
\begin{equation}\label{eq:8}
  \tilde{\mathcal T}_{\lambda}^{+}f(x):= (\lambda
    f_0(x+1) - x f_0(x))\,\mathbb{I}_\N(x),
\end{equation}
which is  none other than the standard
operator for the Poisson distribution. Most authors refer to \eqref{eq:8} as \emph{the}
Stein operator for the Poisson distribution although there are, of course,
many more  operators for this distribution which can be obtained from
\eqref{operator}. One can, for instance,  change the parameterization
of the class $\mathcal{F}(\lambda)$ through
``pre-multiplication'' of the form $f(x) = c(x) f_0(x)$.  See
\cite{GoRe12} for more on this approach. Another way of constructing
Stein operators is by making use of the backward difference, for which  the class
  $\mathcal{F}^{-}(p) =: \mathcal{F}^{-}(\lambda)$ is composed
  of all functions $f:\Z\rightarrow\R$ such that (i) $x\mapsto
  \Delta^{-}(f(x)p_\lambda(x))$ is summable over $\N$ and (ii) 
  $\lim_{x\rightarrow\infty}f(x)p_\lambda(x)=0$. Here no border
  condition   is necessary because $p_\lambda(-1) = 0$. For such $f$ the operator becomes,
  after simplification, 
  \begin{equation*}
    \label{eq:22}
    \mathcal{T}_\lambda^{-}f(x) = \left( f(x) - \frac{x}{\lambda}f(x-1) \right)
 \mathbb{I}_{\N}(x)  \end{equation*}
which is, up to a scaling and a shift, equivalent to the standard
operator \eqref{eq:8}.

Next let $p$ be the  density of $S_n$, the number of
  white balls added to the P\'olya-Eggenberger urn by time $n$,
  with initial state $\alpha\ge1$ white and $\beta\ge1$ black
  balls. We know, e.g. from \cite{GoRe12}, that
  \begin{equation*}
    p(k)=\mathrm{P}\left( S_n=k \right) = \binom{n}{k} \frac{\left(
        \alpha \right)_k \left( \beta \right)_{n-k}}{\left( \alpha+ \beta \right)_n}
  \end{equation*}
  for $k=0,\ldots,n$, with $(x)_0=1$ and otherwise $(x)_k = x(x+1)
  \cdots (x+k-1)$ the rising factorial.  Writing out the classes
  $\mathcal{F}^{\eta}(p)$ and the operators \eqref{operator} in all
  generality for these distributions is of little practical or
  theoretical interest; in particular the resulting objects are hard
  to manipulate (see the discussion in \cite{GoRe12}). It is much more
  informative to directly restrict one's attention to specific
  subclasses. For instance it is easy to see that
  $\mathcal{F}^{+}(p)=:\mathcal{F}^+(\alpha,\beta)$ contains all functions of the form
  $f(x)=xf_0(x)$ with $f_0$ bounded and, for these $f$, the operator
  is of the form
  \begin{equation*}
    \tilde{\mathcal{T}}^{+}_{(\alpha, \beta)}f(x) = \left(
    \frac{(n-x)(\alpha+x)}{\beta+n-x-1}f_0(x+1) - xf_0(x)\right)\mathbb{I}_{[0, n]}(x).
  \end{equation*}
Likewise   $\mathcal{F}^{-}(p)=:\mathcal{F}^{-}(\alpha,\beta)$
 contains all functions of the form 
  $f(x)=(n-x)f_0(x)$ with $f_0$ bounded and, for these $f$, the operator
  is of the form
  \begin{equation*}
    \tilde{\mathcal{T}}^{-}_{(\alpha, \beta)}f(x) = \left(
     (n-x)f_0(x) - f_0(x-1)\frac{x}{\alpha+x-1}(\beta+n-x)\right)\mathbb{I}_{[0, n]}(x).
  \end{equation*}
Of course many variations on the above are imaginable. For
instance one could also choose to consider functions of the form 
 $f(x) = x (\beta+n-x)f_0(x)$; plugging
these into \eqref{operator} yields the operator discussed in
\cite[equation 7]{GoRe12}. 


Thirdly we consider $p$ belonging to the Ord family of distributions,
that is we suppose that there exist $s(x)$ and $\tau(x)$ such that 
\begin{equation*}
  \label{eq:30}
  \frac{p(x+1)}{p(x)} = \frac{s(x) + \tau(x)}{s(x+1)},
\end{equation*}
with $s(a)=0$ (if finite) and $s(x)>0$ for $a<x\le b$. For an
explanation on these notations see \cite[equations (11) and
(12)]{S01}. Writing out the classes $\mathcal{F}^{\eta}(p)$ and the
operators \eqref{operator} in all generality is again of little
practical or theoretical interest.  Note however that
$\mathcal{F}^{+}(p) =: \mathcal{F}^{+}(s,\tau)$ contains all functions
$f:\Z\rightarrow\R$ 
which are of the form $f(x)=f_0(x)s(x)$ with $f_0$ some bounded
function. 
For these $f$, the
operator writes out
$$\tilde{\mathcal T}_{(s, \tau)}^{+}f(x) =
\left((s(x)+\tau(x))f_0(x+1)- s(x) f_0(x)\right)\mathbb{I}_{[a,b]}(x),$$
and we retrieve 
the operator presented in
\cite{S01}.  Similarly for the backward operator we see that 
$\mathcal{F}^{-}(p) =: \mathcal{F}^{-}(s,\tau)$ contains
all functions  $f:\Z\rightarrow\R$ such that (i) $x\mapsto f(x)p(x)$
is bounded over $S_p$ and (ii) $\lim_{x\rightarrow
  b}f(x)p(x)=0$. For these $f$, the operator writes out 
$$\tilde{\mathcal T}_{(s, \tau)}^{-}f(x) = \left( f(x)- \frac{s(x)}{s(x-1)+\tau(x-1)}f(x-1) \right)\mathbb{I}_{[a,b]}(x).$$
There are, of course,   many variations on the approaches presented above. 

Consider next any distribution $p$ on $[0, n]$ satisfying the recurrence 
\begin{equation}\label{eq:larjasstr} 
  a(x)p(x-1) = b(x) p(x) \mbox{ for all } x\in \Z
\end{equation}
with $a(x)$ and $b(x)$ some functions such that $a(x)\neq 0$ for all
$x \in [0, n]$ and $b(0)=0$. Suppose furthermore that $a(n+1) = 0$ (if $n$ is
finite). Then $\mathcal{F}^+(p)$ contains all  functions  of the form
$f(x) = b(x)f_0(x)$ with $f_0$ some bounded function. For these $f$,
the operator writes out 
\begin{equation*}
\mathcal{T}_{(a, b)}^{+}f(x) = \left(a(x+1)f_0(x+1) - f_0(x)b(x)\right)\mathbb{I}_{[0,n]}(x),
\end{equation*}
and we hereby recover \cite[Lemma 2.1]{fulman2012stein}. The specific
distributions studied in Section \ref{sec:other-inequalities} are
obtained by taking  
\begin{equation*}
  a_n(x) = \theta (1-\theta^{-n+x-1}) \mbox{ and } b_n(x) = (\theta^x-1)^2
\end{equation*}
 (distribution of $Q_\theta^n$) and
\begin{equation*}
  a(x) = \theta \mbox{ and } b(x) = (\theta^x-1)^2,
\end{equation*}
 (distribution of $Q_\theta$).

Finally  choose $p$ with support $[0, N]$ for some $N>0$ and represent it
as a Gibbs measure, that is,   write 
\begin{equation*}
  p(x) = \frac{e^{V(x)} \omega^x}{x! \mathcal{Z}}\mathbb{I}_{[0,N]}(x)
\end{equation*}
 with $N$ some positive integer,
  $\omega>0$ fixed, $V$ a function mapping $[0,N]$ to $\R$ and $ V(k)=-\infty$ for $k>N$, and $\mathcal
  Z$ the normalizing constant. This is always possible,
  although there is no unique choice of representation (see
  \cite{ER08}). Then $\mathcal{F}^{\eta}(p) =: \mathcal{F}^{\eta}(V,\omega)$ is
  composed of all functions $f:\Z\rightarrow\R$ which satisfy the
  summability requirements and such that either $f(0)p(0) = 0$ (if
  $\eta=1$) or $f(N)p(N) = 0$ (if $\eta=-1$). In particular, $
  \mathcal{F}^{+}(V,\omega)$
  contains  functions
  of the form 
  $f(x)=xf_0(x)$   with $f_0$ bounded and, for these $f$, the operator
  is of the form
  \begin{equation}
    \label{eq:38}
    \tilde{\mathcal T}^{+}_{(V, \omega)}f(x) =
  \left(e^{V(x+1)-V(x)} \omega f_0(x+1)- xf_0(x)\right)\mathbb{I}_{[0,
    N]}(x);
  \end{equation}
 this corresponds to the Stein operator presented in
  \cite{ER08}. Likewise if $N<\infty$ then $\mathcal{F}^{-}(V,\omega)$
  contains functions of the form $f(x) = (N-x) f_0(x)$ with $f_0$
  bounded and, for these $f$, the operator is of the form
  $$\tilde{\mathcal T}^{-}_{(V, \omega)}f(x) =
  \left(f_0(x)(N-x)- x(N-x+1)\frac{e^{V(x-1)-V(x)}}{\omega}
    f_0(x-1)\right)\mathbb{I}_{[0, N]}(x)$$ and, if $N = \infty$, then
  $f(x) = f_0(x)$ with $f_0$ bounded suffices and the operator is
  equivalent to \eqref{eq:38}. Again a number of other
  parameterizations of the class $\mathcal{F}^{\eta}(V, \omega)$ can
  be considered, each leading to an alternative form of operator.

\section*{Acknowledgments}
\label{sec:acknowledgments}
We thank two anonymous referees and the Associate Editor for their pertinent remarks which have led to
substantial improvement of the paper.

 \bibliographystyle{abbrv}

\bibliography{../../../bibliography/biblio_ys_stein}

\begin{thebibliography}{10}

\bibitem{APP07}
G.~Afendras, N.~Papadatos, and V.~Papathanasiou.
\newblock The discrete {M}ohr and {N}oll inequality with applications to
  variance bounds.
\newblock {\em Sankhy\=a}, 69(2):162--189, 2007.

\bibitem{APP11}
G.~Afendras, N.~Papadatos, and V.~Papathanasiou.
\newblock An extended {S}tein-type covariance identity for the {P}earson family
  with applications to lower variance bounds.
\newblock {\em Bernoulli}, 17(2):507--529, 2011.

\bibitem{BaBaNa03}
K.~Ball, F.~Barthe, and A.~Naor.
\newblock Entropy jumps in the presence of a spectral gap.
\newblock {\em Duke Math. J.}, 119(1):41--63, 2003.

\bibitem{BC05}
A.~D. Barbour and L.~H.~Y. Chen.
\newblock {\em An introduction to Stein's method}, volume~4 of {\em Lect. Notes
  Ser. Inst. Math. Sci. Natl. Univ. Singap.}
\newblock Singapore University Press, Singapore, 2005.

\bibitem{BaHoJa92}
A.~D. Barbour, L.~Holst, and S.~Janson.
\newblock {\em Poisson approximation}, volume~2 of {\em Oxford Studies in
  Probability}.
\newblock The Clarendon Press Oxford University Press, New York, 1992.
\newblock Oxford Science Publications.

\bibitem{BaJoKoMa10}
A.~D. Barbour, O.~Johnson, I.~Kontoyiannis, and M.~Madiman.
\newblock Compound {P}oisson approximation via information functionals.
\newblock {\em Electron. J. Probab.}, 15:1344--1368, 2010.

\bibitem{BA86}
A.~R. Barron.
\newblock Entropy and the central limit theorem.
\newblock {\em Ann. Probab.}, 14(1):336--342, 1986.

\bibitem{Br82}
L.~D. Brown.
\newblock A proof of the central limit theorem motivated by the
  {C}ram{\'e}r-{R}ao inequality.
\newblock In {\em Statistics and probability: essays in honor of {C}. {R}.
  {R}ao}, pages 141--148. North-Holland, Amsterdam, 1982.

\bibitem{carlen1991entropy}
E.~Carlen and A.~Soffer.
\newblock Entropy production by block variable summation and central limit
  theorems.
\newblock {\em Commun. Math. Phys.}, 140(2):339--371, 1991.

\bibitem{ChGoSh11}
L.~H.~Y. Chen, L.~Goldstein, and Q.-M. Shao.
\newblock {\em Normal approximation by {S}tein's method}.
\newblock Probability and its Applications (New York). Springer, Heidelberg,
  2011.

\bibitem{CoTh06}
T.~Cover and J.~Thomas.
\newblock {\em Elements of Information Theory}, volume Second Edition.
\newblock Wiley \& Sons, New York, 2006.

\bibitem{ER08}
P.~Eichelsbacher and G.~Reinert.
\newblock Stein's method for discrete gibbs measures.
\newblock {\em Ann. Appl. Probab.}, 18:1588--1618, 2008.

\bibitem{E05}
T.~Erhardsson.
\newblock Stein's method for poisson and compound poisson approximation.
\newblock In {\em An introduction to Stein's method}, 2005.

\bibitem{fulman2012stein}
J.~Fulman and L.~Goldstein.
\newblock Stein's method and the rank distribution of random matrices over
  finite fields.
\newblock {\em arXiv preprint arXiv:1211.0504}, 2012.

\bibitem{GS02}
A.~L. Gibbs and F.~E. Su.
\newblock On choosing and bounding probability metrics.
\newblock {\em International Statistical Review / Revue Internationale de
  Statistique}, 70(3):pp. 419--435, 2002.

\bibitem{GoRe12}
L.~Goldstein and G.~Reinert.
\newblock Stein's method and the beta distribution.
\newblock Preprint, arxiv:1207.1460, 2012.

\bibitem{hillion2011natural}
E.~Hillion, O.~Johnson, and Y.~Yu.
\newblock A natural derivative on [0, n] and a binomial poincar$\backslash$'e
  inequality.
\newblock {\em arXiv preprint arXiv:1107.0127}, 2011.

\bibitem{Ho04}
S.~Holmes.
\newblock Stein's method for birth and death chains.
\newblock In {\em Stein's method: expository lectures and applications},
  volume~46 of {\em IMS Lecture Notes Monogr. Ser.}, pages 45--67. Inst. Math.
  Statist., Beachwood, OH, 2004.

\bibitem{Jo04}
O.~Johnson.
\newblock {\em Information theory and the central limit theorem}.
\newblock Imperial College Press, London, 2004.

\bibitem{MR2128239}
O.~Johnson and A.~Barron.
\newblock Fisher information inequalities and the central limit theorem.
\newblock {\em Probab. Theory Related Fields}, 129(3):391--409, 2004.

\bibitem{JM87}
I.~Johnstone and B.~MacGibbon.
\newblock Une mesure d'information caract\'erisant la loi de poisson.
\newblock In {\em S{\'e}minaire de probabilit{\'e}s}, volume XXI, pages
  563--573. Springer, 1987.

\bibitem{KoHaJo05}
I.~Kontoyiannis, P.~Harremo{{\"e}}s, and O.~Johnson.
\newblock Entropy and the law of small numbers.
\newblock {\em IEEE Trans. Info. Theory}, 51:466--472, 2005.

\bibitem{KullIEEE}
S.~Kullback.
\newblock A lower bound for discrimination information in terms of variation.
\newblock {\em IEEE Trans. Info. Theory}, 4, 1967.

\bibitem{LS12a}
C.~Ley and Y.~Swan.
\newblock Stein's density approach and information inequalities.
\newblock {\em Electron. Comm. Probab.}, 18(7):1--14, 2013.

\bibitem{NP11}
I.~Nourdin and G.~Peccati.
\newblock {\em Normal approximations with Malliavin calculus : from Stein's
  method to universality}.
\newblock Cambridge Tracts in Mathematics. Cambridge University Press, 2012.

\bibitem{nourdin2013entropy}
I.~Nourdin, G.~Peccati, and Y.~Swan.
\newblock Entropy and the fourth moment phenomenon.
\newblock {\em arXiv preprint arXiv:1304.1255}, 2013.

\bibitem{sason2012entropy}
I.~Sason.
\newblock Entropy bounds for discrete random variables via coupling.
\newblock Preprint, arXiv:1209.5259, 2012.

\bibitem{sason2012information}
I.~Sason.
\newblock An information-theoretic perspective of the poisson approximation via
  the chen-stein method.
\newblock {\em arXiv preprint arXiv:1206.6811}, 2012.

\bibitem{sason2012entropybern}
I.~Sason.
\newblock On the entropy of sums of bernoulli random variables via the
  chen-stein method.
\newblock In {\em Information Theory Workshop (ITW), 2012 IEEE}, pages
  542--546. IEEE, 2012.

\bibitem{sason2013improved}
I.~Sason.
\newblock Improved lower bounds on the total variation distance and relative
  entropy for the poisson approximation.
\newblock {\em arXiv preprint arXiv:1301.7504}, 2013.

\bibitem{S01}
W.~Schoutens.
\newblock Orthogonal polynomials in {S}tein's method.
\newblock {\em J. Math. Anal. Appl.}, 253(2):515--531, 2001.

\bibitem{Sh75}
R.~Shimizu.
\newblock On {F}isher's amount of information for location family.
\newblock In {\em A Modern Course on Statistical Distributions in Scientific
  Work}, pages 305--312. Springer, 1975.

\end{thebibliography}

\end{document}